\documentclass[11pt]{amsart}
\usepackage{amssymb}
\usepackage{amsmath}
\usepackage{youngtab}

\begin{document}  
\newtheorem{defn0}{Definition}[section]
\newtheorem{prop0}[defn0]{Proposition}
\newtheorem{thm0}[defn0]{Theorem}
\newtheorem{lemma0}[defn0]{Lemma}
\newtheorem{coro0}[defn0]{Corollary}
\newtheorem{exa}[defn0]{Example}
\newtheorem{exe}[defn0]{Exercise}
\newtheorem{rem0}[defn0]{Remark}
\def\rig#1{\smash{ \mathop{\longrightarrow}\limits^{#1}}}
\def\swar#1{\swarrow\rlap{$\vcenter{\hbox{$\scriptstyle#1$}}$}}
\def\lswar#1{\swarrow\llap{$\vcenter{\hbox{$\scriptstyle#1$}}$}}
\def\sear#1{\searrow\rlap{$\vcenter{\hbox{$\scriptstyle#1$}}$}}
\def\lsear#1{\searrow\llap{$\vcenter{\hbox{$\scriptstyle#1$}}$}}
\def\near#1{\nearrow\rlap{$\vcenter{\hbox{$\scriptstyle#1$}}$}}
\def\dow#1{\Big\downarrow\rlap{$\vcenter{\hbox{$\scriptstyle#1$}}$}}
\def\ldow#1{\Big\downarrow\llap{$\vcenter{\hbox{$\scriptstyle#1$}}$}}
\def\up#1{\Big\uparrow\rlap{$\vcenter{\hbox{$\scriptstyle#1$}}$}}
\def\lef#1{\smash{ \mathop{\longleftarrow}\limits^{#1}}}
\newcommand{\defref}[1]{Def.~\ref{#1}}
\newcommand{\rank}{\operatorname{rank}}
\newcommand{\propref}[1]{Prop.~\ref{#1}}
\newcommand{\thmref}[1]{Thm.~\ref{#1}}
\newcommand{\lemref}[1]{Lemma~\ref{#1}}
\newcommand{\corref}[1]{Cor.~\ref{#1}}
\newcommand{\exref}[1]{Example~\ref{#1}}
\newcommand{\secref}[1]{Section~\ref{#1}}
\newcommand{\qedd}{\hfill\framebox[2mm]{\ }\medskip}
\def\P#1{{\bf P}^#1}
\def\I{{\cal I}}
\def\O{{\mathcal O}}
\def\C{{\mathbb C}}
\def\R{{\mathbb R}}
\def\Z{{\mathbb Z}}
\def\proof{{\it Proof:\  \ }}
\def\M{{\cal M}}
\def\F{{\cal F}}

\title{An invariant regarding Waring's problem\\ for cubic polynomials }

\author{ Giorgio Ottaviani}

\dedicatory{to the memory of Michael Schneider, ten years after}


\begin{abstract}

\noindent 
We compute the equation of the $7$-secant variety to the Veronese variety
$(\P 4,\O(3))$, its degree is $15$. This is the last missing invariant in the Alexander-Hirschowitz classification.
 It gives the condition to express a homogeneous cubic polynomial in $5$ variables as the sum of $7$ cubes
 (Waring problem). The interesting side in the construction is that it comes from the determinant of
 a matrix of order $45$ with linear entries, which is a cube. The same technique allows to express the classical Aronhold invariant
of plane cubics as a pfaffian.   

\end{abstract}

\maketitle

\section{Introduction}

We work over an algebraically closed field $K$ of characteristic zero.
The Veronese variety,
given by $\P n$ embedded with the linear system $|\O(d)|$,
lives in $\P N$ where $N=\binom{n+d}{d}-1$. It parametrizes the
homogeneous polynomials $f$ of degree $d$ in $n+1$ variables
which are the power of a linear form $g$, that is $f=g^d$.

Let $\sigma_s(\P n,\O(d))$ be the $s$-secant variety of the Veronese variety, that is
the Zariski closure of the variety of polynomials $f$ which are the sum
of the powers of $s$ linear forms $g_i$, i.e. $f=\sum_{i=1}^sg_i^d$.
In particular $\sigma_1(\P n,\O(d))$ is the Veronese variety itself and
$\sigma_2(\P n,\O(d))$ is the usual secant variety. For generalities
about the Waring's problem for polynomials see \cite{IK} or \cite{RS}.

Our starting point is the theorem of Alexander and Hirschowitz (see \cite{AH} or \cite{BO}
for a survey, including a self-contained proof)
which states that the codimension of  $\sigma_s(\P n,\O(d))\subseteq\P N$
is the expected one, that is $\max\{N+1-(n+1)s,0\}$, with the only exceptions

\begin{itemize}
\item{(i)} $\sigma_k(\P n,\O(2)),\ 2\le k\le n$
\item{(ii)} $\sigma_{\frac{1}{2}n(n+3)}(\P n,\O(4)),\ n=2,3,4$ 
\item{(iii)} $\sigma_7(\P 4,\O(3))$ 
\end{itemize}
The case (i)  corresponds to the matrices of rank $\le k$ in the variety of symmetric matrices
of order $n+1$.  In the  cases (ii) and (iii) the expected codimension is zero, while the codimension is one.
Hence the equation of the hypersurface $\sigma_s(\P n,\O(d))$ in these cases is an interesting $SL(n+1)$-invariant.
 In the cases  (ii) it is the catalecticant invariant, that was computed by Clebsch in the XIX century,
 its degree is $\binom{n+2}{2}$.
 
The main result of this paper is the computation of the equation of $\sigma_7(\P 4,\O(3))$.
This was left as an open problem in \cite[chap.2, rem. 2.4]{IK} .

We consider a vector space $V$. For any nonincreasing sequence of positive integers 
$\alpha=(\alpha_1,\alpha_2,\ldots)$
it is defined the Schur module $\Gamma^{\alpha}V$, which is an irreducible $SL(V)$-module (see \cite{FH}). For $\alpha=(p)$ we get the $p$-th symmetric power of $V$
and for $\alpha=(1,\ldots ,1)$ ($p$ times) we get the $p$-th alternating power of $V$. 
The module $\Gamma^{\alpha}V$
is visualized as a Young diagram containing $\alpha_i$ boxes in the $i$-th row. In particular 
if $\dim V=5$ then
$\Gamma^{2,2,1,1}V$ and its dual $\Gamma^{2,1,1}V$ have both dimension $45$. 

Our main result is the following

\begin{thm0}\label{main}
Let $V$ be a vector space of dimension $5$. For any $\phi\in S^3V$, let
$B_{\phi}\colon\Gamma^{2,2,1,1}V\to\Gamma^{2,1,1}V$ be the $SL(V)$-invariant contraction operator.
Then there is an irreducible homogeneous polynomial  $P$ of degree $15$  on $S^3V$
such that $$2P(\phi)^3=\det B_{\phi}$$
The polynomial $P$ is the equation of $\sigma_7(\P {{}}(V),\O(3))$.
\end{thm0}

The coefficient $2$ is needed because we want the invariant polynomials to be defined over the rational numbers.
The picture in terms of Young diagrams is
$$\yng(2,2,1,1)\otimes\young(***)\quad\to\quad\young(\ \ *,\ \ ,\ *,\ ,*)\simeq\ \yng(2,1,1)$$

This picture means that $\Gamma^{2,1,1}V$ is a direct summand
of the tensor product $\Gamma^{2,2,1,1}V\otimes S^3V$, according to the Littlewood-Richardson rule (\cite{FH}).

The polynomial $P$ gives the necessary condition  to express a cubic homogeneous polynomial
in five variables as a sum of seven cubes. We prove in \lemref{im2} that 
if $\phi$ is decomposable then $\textrm{rk}(B_{\phi})=6$.
 The geometrical explanation that $\sigma_7(\P 4,\O(3))$ is an exceptional case is related to the fact
 that given seven points in $\P 4$ there is a unique rational normal curve through them,
and it was discovered independently by Richmond and Palatini in 1902,
see \cite{CH} for a modern reference.
Our approach gives a different (algebraic) proof of the fact that
$\sigma_7(\P 4,\O(3))$ is an exceptional case. Another argument, by using syzygies, is in \cite{RS}.
 B. Reichstein found in 
\cite{Re} an algorithm to check when a cubic homogeneous polynomial
in five variables is the sum of seven cubes, see the Remark \ref{reich}.

The resulting table of the Alexander-Hirschowitz classification  is the following

\[\begin{array}{|l|l|l|l|}
\hline
&exp. codim&codim&equation\\
\hline 
\sigma_k(\P n,\O(2))\quad 2\le k\le n &\max(\frac{(n+1)(n+2-2k)}{2},0)&\binom{n-k+2}{2}&
(k+1)-\textrm{minors }
 \\
\hline 
\sigma_{\frac{1}{2}n(n+3)}(\P n,\O(4))\quad n=2,3,4 &0&1&\textrm{catalecticant inv.} \\
\hline 
\sigma_7(\P 4,\O(3))&0&1&\textrm{see \thmref{main}} \\
\hline
\end{array}\]

The degree of $\sigma_k(\P n,\O(2))$ was computed by C. Segre, it is equal to
  $\prod_{i=0}^{n-k}\frac{\binom{n+1+i}{n+1-k-i}}{\binom{2i+1}{i}}$.
We will use in the proof of \thmref{main}
the fact that $\sigma_{k-1}(\P n,\O(2))$ is the singular locus of
$\sigma_k(\P n,\O(2))$ for $k\le n$.

A general cubic polynomial
in five variables can be expressed as a sum of eight cubes in $\infty^5$ ways,
parametrized by a Fano $5$-fold of index one (see \cite{RS}).
A  cubic polynomial
in five variables  which can be expressed as a sum of seven cubes was called degenerate
in \cite{RS},
hence what we have found is the locus of degenerate cubics. A degenerate cubic
in five variables can be expressed as a sum of seven cubes in $\infty^1$ ways,
parametrized by ${\bf P}^1$ (see \cite{RS} 4.2).

 To explain our technique, we consider the Aronhold
 invariant of plane cubics.

The Aronhold invariant is the degree $4$ equation of $\sigma_3(\P 2,\O(3))$,
which can be seen as the $SL(3)$-orbit of the Fermat cubic $x_0^3+x_1^3+x_2^3$ (sum of three cubes),
see \cite[Prop. 4.4.7]{St}  or \cite[(5.13.1)]{DK} .
  
Let $W$ be a vector space of dimension $3$. In particular 
$\Gamma^{2,1}W=\textrm{ad\ }W$ is self-dual and it has dimension $8$. We get
\begin{thm0}\label{aronhold}
For any $\phi\in S^3W$, let
$A_{\phi}\colon\Gamma^{2,1}W\to\Gamma^{2,1}W$ be the $SL(V)$-invariant contraction operator.
Then  $A_{\phi}$ is skew-symmetric and the pfaffian $\textrm{Pf\ } A_{\phi}$
 is the equation of $\sigma_3(\P {{}}(W),\O(3))$,
i.e. it is the Aronhold invariant.
\end{thm0}

The corresponding picture is
$$\yng(2,1)\otimes\young(***)\quad\to\quad\young(\ \ *,\ *,*)\simeq\ \yng(2,1)$$

The Aronhold invariant gives the necessary condition  to express a cubic homogeneous polynomial
in three variables as a sum of three cubes. The explicit expression of the Aronhold invariant
 is known since the XIX century, but we have not found in the literature its representation as a pfaffian.
In the remark \ref{scorza} we apply this representation to the Scorza map between plane quartics.

In section 2 we give the proof of \thmref{aronhold}. This is introductory to \thmref{main},
which is proved in section 3. In section 4 we review, for completeness, some known facts about the 
catalecticant invariant of quartic hypersurfaces.

  We are indebted to S. Sullivant, for his beautiful lectures
at Nordfjordeid in 2006 about \cite{SS}, where a representation of the Aronhold invariant is found
with combinatorial techniques.

 \section{The Aronhold invariant as a pfaffian}
 
Let $e_0, e_1, e_2$ be a basis of $W$ and fix the orientation
$\wedge^3W\simeq K$ given by $e_0\wedge e_1\wedge e_2$. We have $\textrm{End}~W=\textrm{ad}~W\oplus K$.
 The $SL(W)$-module $\textrm{ad}~W=\Gamma^{2,1}(W)$
 consists of the subspace of endomorphisms of $W$ with zero trace.
  We may interpret the contraction 
  $$A_{\phi}\colon\Gamma^{2,1}W\to\Gamma^{2,1}W$$  as the restriction of a linear map 
  $A'_{\phi}\colon \textrm{End}~W\to \textrm{End}~W$, which is defined 
  for $\phi=e_{i_1}e_{i_2}e_{i_3}$ as
 $$A'_{e_{i_1}e_{i_2}e_{i_3}}(M)(w)=
\sum_{\sigma}(M(e_{i_{\sigma(1)}})\wedge e_{i_{\sigma(2)}}\wedge w)e_{i_{\sigma(3)}}$$
 where $M\in \textrm{End}~W$, $w\in W$ and $\sigma$ covers the symmetric group $\Sigma_3$.
 
 Then $A'_{\phi}$ is defined for a general $\phi$ by linearity,
 and it follows from the definition that it is $SL(V)$-invariant. 
 
 The Killing scalar product
 on $\textrm{End} W$ is defined by $tr(M\cdot N)$.
 
 \begin{lemma0}\label{skew}
 (i) $\textrm{Im}(A'_{\phi})\subseteq \textrm{ad}~W\qquad K\subseteq\textrm{Ker}(A'_{\phi})$

(ii)  $A'_{\phi}$ is skew-symmetric.
 \end{lemma0}
 
 {\it Proof} (i) follows from
  $$ tr\left[A_{e_{i_1}e_{i_2}e_{i_3}}(M)\right]=\sum_s A_{e_{i_1}e_{i_2}e_{i_3}}(M)(e_s)e_s^{\vee}=
 \sum_{\sigma}(M(e_{i_{\sigma(1)}})\wedge e_{i_{\sigma(2)}}\wedge e_{i_{\sigma(3)}})=0$$
The second inclusion is evident.
  To prove (ii), we have to check that
 $$tr(A_{\phi}(M)\cdot N)=-tr(A_{\phi}(N)\cdot M)$$
 for $M, N\in \textrm{End}~W$. Indeed let $\phi=e_{i_1}e_{i_2}e_{i_3}$. We get 
  $$tr(A_{e_{i_1}e_{i_2}e_{i_3}}(M)\cdot N)=\sum_s A_{e_{i_1}e_{i_2}e_{i_3}}(M)(N(e_s))e_s^{\vee}=
  \sum_{\sigma}M(e_{i_{\sigma(1)}})\wedge e_{i_{\sigma(2)}}\wedge N(e_{i_{\sigma(3)}})$$
  which is alternating in $M$ and $N$, where we denoted by $e_i^{\vee}$ the dual basis.\qedd

It follows from \lemref{skew} that the restriction
$${A'_{\phi}}_{|\textrm{ad}~W}\colon\textrm{ad}~W\to\textrm{ad}~W$$ coincides, up to scalar multiple,
with the contraction operator $A_{\phi}$ of \thmref{aronhold} and it is skew-symmetric.
  
\begin{lemma0}\label{im}
Let
 $\phi=w^3$ with $w\in W$. Then $\textrm{rk}A_{\phi}=2$. More precisely 
 $$\textrm{Im}A_{w^3}=\{M\in \textrm{ad}~W| \textrm{Im}~M\subseteq <w>\}$$
 $$\textrm{Ker}A_{w^3}=\{M\in \textrm{ad}~W| w\textrm{\ is an eigenvector of\ }M\}$$
\end{lemma0}

 {\it Proof} The statement follows from the equality
 $$A_{w^3}(M)(v)=6(M(w)\wedge w\wedge v)w$$
 As an example, note that $\textrm{Im} A_{e_0^3}=<e_0\otimes e_1^{\vee}, e_0\otimes e_2^{\vee}>$
and $\textrm{Ker} A_{e_0^3}$ is spanned by all the basis monomials, with the exception
of $e_0^{\vee}\otimes e_1$ and $e_0^{\vee}\otimes e_2$.
 Due to the $SL(W)$-invariance, this example proves the general case. \qedd
 
 {\it Proof of \thmref{aronhold}}
 
 Let $\phi\in\sigma_3(\P{{}}(W),\O(3))$. By the definition of higher secant variety, $\phi$ is
 in the closure of elements
 which can be written as $\phi_1+\phi_2+\phi_3$ with $\phi_i\in (\P{{}}(W),\O(3))$.
 From \lemref{im}
 it follows that  $$\textrm{rk\ }A_{\phi}\le \textrm{rk\ }A_{\sum_{i=1}^3\phi_i}=
 \textrm{rk\ }\sum_{i=1}^3 A_{\phi_i}\le \sum_{i=1}^3\textrm{rk\ }A_{\phi_i}= 2\cdot 3 = 6$$
 Hence $Pf(A_{\phi})$ has to vanish on $\sigma_3(\P{{}}(W),\O(3))$.
 
 Write a cubic polynomial as

$$\phi = v_{000}x_0^3+3v_{001}x_0^2x_1+3v_{002}x_0^2x_2+3v_{011}x_0x_1^2+6v_{012}x_0x_1x_2+3v_{022}x_0x_2^2+$$
$$+v_{111}x_1^3+3v_{112}x_1^2x_2+3v_{122}x_1x_2^2+v_{222}x_2^3$$

We order the monomial basis of $\wedge^2W\otimes W$ with the lexicographical order in the following way:

$(w_0\wedge w_1)w_0, (w_0\wedge w_1)w_1, (w_0\wedge w_1)w_2, (w_0\wedge w_2)w_0, (w_0\wedge w_2)w_1, (w_0\wedge w_2)w_2,$

$ (w_1\wedge w_2)w_0, (w_1\wedge w_2)w_1, (w_1\wedge w_2)w_2$

Call $M_i$ for $i=1,\ldots ,9$ this basis.
The matrix of $A'_{\phi}$,  with respect to this basis, has at the entry $(i,j)$
the value $A'_{\phi}(M_j)(M_i)$ and it is the following
$$\left[\begin{array}{rrrrrrrrr}
0&v_{222}&-v_{122}&0&-v_{122}&v_{112}&0&v_{022}&-v_{012}\\
-v_{222}&0&v_{022}&v_{122}&0&-v_{012}&-v_{022}&0&v_{002}\\
v_{122}&-v_{022}&0&-v_{112}&v_{012}&0&v_{012}&-v_{002}&0\\
0&-v_{122}&v_{112}&0&v_{112}&-v_{111}&0&-v_{012}&v_{011}\\
v_{122}&0&-v_{012}&-v_{112}&0&v_{011}&v_{012}&0&-v_{001}\\
-v_{112}&v_{012}&0&v_{111}&-v_{011}&0&-v_{011}&v_{001}&0\\
0&v_{022}&-v_{012}&0&-v_{012}&v_{011}&0&v_{002}&-v_{001}\\
-v_{022}&0&v_{002}&v_{012}&0&-v_{001}&-v_{002}&0&v_{000}\\
v_{012}&-v_{002}&0&-v_{011}&v_{001}&0&v_{001}&-v_{000}&0\\
\end{array}
\right]$$

Deleting one of the columns corresponding to $(w_0\wedge w_1)w_2$, $(w_0\wedge w_2)w_1$
or $(w_1\wedge w_2)w_0$ (respectively the $3$rd, the $5$th and the $7$th, indeed their alternating
 sum gives the trace), and the corresponding row, we get a skew-symmetric matrix of order $8$
which is the matrix of $A_{\phi}$.  To conclude the proof,
 it is enough to check that the pfaffian is nonzero.
This can be easily checked on the point corresponding
 to $\phi=x_0x_1x_2$, that is when $v_{012}=1$ and all the other coordinates are equal to zero.
 This means that any triangle is not in the closure of the Fermat curve.
we conclude that $Pf(A_{\phi})$ is the Aronhold invariant. We verified that it
coincides, up to a constant, with the expression given in
 \cite[Prop. 4.4.7]{St}  or in \cite[(5.13.1)]{DK} . 
 
\qedd

 The vanishing of the Aronhold invariant gives the necessary and sufficient condition 
to express a cubic polynomial in three variables as the sum of three cubes.

 {\bf Remark} $A'_{\phi}$ can be thought as a map
 $$A'_{\phi}\colon\wedge^2W\otimes W\to \wedge^2W^{\vee}\otimes W^{\vee}$$
 For $\phi=w^3$ we have the formula
 $$A'_{\phi}(\omega\otimes v)(\omega'\otimes v')=
 (\omega\wedge w)\otimes(v\wedge w\wedge v')\otimes(\omega'\wedge w)$$ 
This is important for the understanding of the next section.

 {\bf Remark} We have the decomposition
$$\wedge^2(\Gamma^{2,1}W)=S^3W\oplus\Gamma^{2,2,2}W\oplus\textrm{ad\ }W$$

and it is a nice exercise to show the behaviour of the three summands.
For the first one
$$S^3W\cap\{M\in \wedge^2(\Gamma^{2,1}W)| \textrm{rk\ }(M)\le 2k\}$$
is the cone over
$\sigma_k\left({\bf P}(W),\O(3)\right)$,
so that we have found the explicit equations for all the higher secant varieties
to  $\left({\bf P}(W),\O(3)\right)$.  The secant variety
$\sigma_2\left({\bf P}(W),\O(3)\right)$ is the closure of the orbit of plane cubics consisting of three concurrent 
lines, and its equations are the $6\times 6$ subpfaffians of $A_{\phi}$. It has degree $15$.
There is a dual description for $\Gamma^{2,2,2}W$.

For the third summand, we have that 
$$\textrm{ad\ }W\subseteq \{M\in \wedge^2(\Gamma^{2,1}W)| \textrm{rk\ }(M)\le 6\}$$
Indeed any $M\in \textrm{ad\ }W$ induces the skew-symmetric morphism
$$[M,-]$$
whose kernel contains  $M$.
Moreover
$$\textrm{ad\ }W\cap \{M\in \wedge^2(\Gamma^{2,1}W)| \textrm{rk\ }(M)\le 4\}$$
is the $5$-dimensional affine cone consisting of endomorphisms $M\in \textrm{ad\ }W$ such that
their minimal polynomial has degree $\le 2$.

\begin{rem0}\label{scorza}  We recall from \cite{DK} the definition of the Scorza map.
Let $A$ be the Aronhold invariant. For any plane quartic $F$ and any point $x\in{\bf P}(W)$
we consider the polar cubic  $P_x(F)$. Then $A(P_x(F))$ is a quartic in the variable $x$
which we denote by $S(F)$. The rational map $S\colon{\bf P}(S^4W)\dashrightarrow {\bf P}(S^4W)$
 is called the Scorza map. Our description of the Aronhold invariant shows
that $S(F)$ is defined as the degeneracy locus of a skew-symmetric morphism  on ${\bf P}(W)$
$$\O (-2)^8\rig{f} \O(-1)^8$$
It is easy to check  (see \cite{Be}) that $\textrm{Coker\ }f=E$ is a rank two vector bundle
over $S(F)$ such that $c_1(E)=K_{S(F)}$. Likely from $E$ it is possible to recover the eben theta-characteristic 
$\theta$ on $S(F)$ defined in \cite[(7.7)]{DK} . The natural guess is that
$$h^0(E\otimes(-\theta))>0$$ for a unique even $\theta$, but we do not know if this is true.
\end{rem0}

\section{The invariant for cubic polynomials in five variables}

Let now $e_0,\ldots ,e_4$ be a basis of $V$, no confusion will arise with the notations of
the previous section.  We fix the orientation
$\wedge^5 V\simeq K$ given by $e_0\wedge e_1\wedge e_2\wedge e_3\wedge e_4$.
We construct, for $\phi\in S^3V$, the contraction operator
$$B'_{\phi}\colon\wedge^4V\otimes\wedge^2V\to \wedge^4V^{\vee}\otimes\wedge^2V^{\vee}\simeq \wedge^3V\otimes V$$
For a decomposable $\phi=e_{i_1}e_{i_2}e_{i_3}$, the definition is
 $$B'_{\phi}(v_a\wedge v_b\wedge v_c\wedge v_d)\otimes(v_e\wedge v_f)=
 \sum_{\sigma} \left(v_a\wedge v_b\wedge v_c\wedge v_d\wedge e_{i_{\sigma(1)}}\right)\otimes
 \left(v_e\wedge v_f\wedge e_{i_{\sigma(2)}}\right)\otimes e_{i_{\sigma(3)}}$$
where $\sigma$ covers the symmetric group $\Sigma_3$
 and we extend this definition, to a general $\phi$, by linearity.
 
We may interpret $B'_{\phi}$ as a morphism
$$B'_{\phi}\colon Hom(V,\wedge^2V)\to Hom(\wedge^2V,V)$$

If $\phi=e_{i_1}e_{i_2}e_{i_3}$ and $M\in Hom(V,\wedge^2V)$ we have

$$B'_{e_{i_1}e_{i_2}e_{i_3}}(M)(v_1\wedge v_2)=\sum_{\sigma}(M(e_{i_{\sigma(1)}})\wedge e_{i_{\sigma(2)}}\wedge v_1\wedge v_2)e_{i_{\sigma(3)}}$$
 
We have a $SL(V)$-decomposition
$$\wedge^4V\otimes\wedge^2V= \Gamma^{2,2,1,1}V\oplus V$$
Consider the contraction
$c\colon\wedge^4V\otimes\wedge^2V\to V$
defined by
$$c\left(\omega\otimes (v_i\wedge v_j)\right) = (\omega\wedge v_i) v_j-(\omega\wedge v_j) v_i$$
Then the subspace $\Gamma^{2,2,1,1}V$ can be identified with
$\{M\in \wedge^4V\otimes\wedge^2V| c(M)=0\}$
or with
$$\{M\in Hom(V,\wedge^2V)| \sum e_i^{\vee}M(e_i)=0\}$$
 The subspace $V\subset Hom(V,\wedge^2V)$ can be identified with $\{v\wedge - | v\in V\}$.
At the same time we have a $SL(V)$-decomposition
$$V\otimes\wedge^3 V= \Gamma^{2,1,1}V\oplus \wedge^4 V$$
and the obvious contraction $d\colon V\otimes\wedge^3V\to\wedge^4V$.
The subspace $\Gamma^{2,1,1}V$ can be identified with
$$\{N\in V\otimes \wedge^3V| d(N)=0\}$$
 
 \begin{lemma0}\label{sym}
 (i) $\textrm{Im}(B'_{\phi})\subseteq \Gamma^{2,1,1}V\qquad V\subseteq\textrm{Ker}(B'_{\phi}) $
 
 (ii) $B'_{\phi}$ is symmetric.
 \end{lemma0}
 
 {\it Proof} The statement (i) follows from the formula
 
 $$d\left(B'_{e_{i_1}e_{i_2}e_{i_3}}(v_a\wedge v_b\wedge v_c\wedge v_d)\otimes(v_e\wedge v_f)\right) =$$
 $$=\sum_{\sigma} \left(v_a\wedge v_b\wedge v_c\wedge v_d\wedge e_{i_{\sigma(1)}}\right)\otimes
 \left(v_e\wedge v_f\wedge e_{i_{\sigma(2)}}\wedge e_{i_{\sigma(3)}}\right)=0$$
  In order to prove the second inclusion, for any $v\in V$ consider the induced morphism $M_v(w)=v\wedge w$.
 We get
$$B'_{e_{i_1}e_{i_2}e_{i_3}}(M_v)(v_1\wedge v_2)=\sum_{\sigma}
\left(v\wedge e_{i_{\sigma(1)}}\wedge e_{i_{\sigma(2)}}\wedge v_1\wedge v_2\right) e_{i_{\sigma(3)}}=0$$

 In order to prove (ii) we may assume
 $\phi=v^3$.
 
 We need to prove that
 $$B'_{v^3}(\omega\otimes \xi)(\omega'\otimes \xi')= B'_{v^3}(\omega'\otimes \xi')(\omega\otimes \xi)$$
  for every $\omega, \omega'\in\wedge^4V$ and $\xi, \xi'\in\wedge^2V$.
 Indeed
 $$B'_{v^3}(\omega\otimes \xi)(\omega'\otimes \xi')=
 (\omega\wedge v)\otimes(\xi\wedge v\wedge \xi')\otimes(v\wedge \omega')$$
 which is symmetric in the pair $(\omega,\xi)$.\qedd

It follows from \lemref{sym} that the restriction
  ${B'_{\phi}}_{|\Gamma^{2,2,1,1}}\colon\Gamma^{2,2,1,1}\to\Gamma^{2,1,1}V$
coincides, up to scalar multiple, with the contraction $B_{\phi}$ of the \thmref{main}
and it is symmetric.
Note that $$\textrm{Ker}(B_{\phi})=\textrm{Ker}(B'_{\phi})/V\qquad \textrm{Im}(B_{\phi})= \textrm{Im}(B'_{\phi})$$

 \begin{lemma0}\label{im2}
Let
 $\phi=v^3$ with $v\in V$. Then $\textrm{rk\ }B_{\phi}=6$. More precisely 
 $$\textrm{Im}B_{v^3}=\{N\in Hom(\wedge^2 V,V) |  \sum e_i^{\vee}N(e_i\wedge v)=0\quad\forall v\in V,\quad Im(N)\subseteq <v>\}$$
 $$\textrm{Ker}B_{v^3}=\{M\in Hom(V,\wedge^2 V) | \sum e_i^{\vee}M(e_i)=0,\quad M(v)\subseteq v\wedge V\}$$
 \end{lemma0}

 {\it Proof\ }  The statement follows from the equality
 $$B_{v^3}(M)(v_1\wedge v_2)=6\left(M(v)\wedge v\wedge v_1\wedge v_2\right)v$$ 
 
 As an example, a basis of $\textrm{Im}B_{e_0^3}$ is given by $e_0\otimes (e_i^{\vee}\wedge e_j^{\vee})$
 for $1\le i<j\le 4$ and a basis of $\textrm{Ker}B_{e_0^3}$
is given by all the basis monomials with the exceptions of $e_0^{\vee}\otimes (e_i\wedge e_j)$
for $1\le i<j\le 4$. Due to the $SL(V)$-invariance,
this example proves the general case. \qedd
 
We write $\phi\in S^3V$ as
$\phi=v_{000}x_0^3+3v_{001}x_0^2x_1+\ldots +v_{444}x_4^3$ 
 
 \begin{lemma0}\label{monomial}
 Every $SL(V)$-invariant homogeneous polynomial of degree $15$ on $S^3V$ which contains the monomial
 $$v_{000}^2v_{012}^3v_{111}v_{223}^3v_{334}^3v_{144}^3$$ is irreducible.
 \end{lemma0}
 
  {\it Proof}  Let $t_0,\ldots ,t_4$ be the canonical basis
  of ${\Z}^5$. We denote by
  $t_i+t_j+t_k$ the weight of the monomial $v_{ijk}$,  according to \cite{St}. For example the weight of
  $v_{000}$ is $(3,0,0,0,0)$ . We denote
the first component of the weight as the $x_0$-weight, the second component as the $x_1$-weight,
and so on. We recall that
  every $SL(V)$-invariant polynomial is isobaric, precisely
  every monomial of a $SL(V)$-invariant polynomial of degree $5k$
  has weight $(3k,3k,3k,3k,3k)$ (see \cite[(4.4.14)]{St} ), this follows from the invariance with respect
to the diagonal torus. We claim that 
  there is no isobaric monomial of weight $(6,6,6,6,6)$ and degree 10 with variables among
 $v_{000}, v_{012}, v_{111}, v_{223}, v_{334}, v_{144}$.
 We divide into the following cases,
 by looking at the possibilities for the $x_0$-weight:
 \begin{enumerate}
 \item[i)] The monomial contains $v_{000}^2$ and does not contain $v_{012}$.
By looking at  the
 $x_2$-weight, the monomial has to contain $v_{223}^3$, which gives contribution $3$ to the $x_3$-weight.
 This gives a contradiction ,
 because from $v_{334}$ the possible values for the $x_3$-weight are even, and we never make $6$.
 \item[ii)] The monomial contains $v_{000}v_{012}^3$ and not higher powers. This monomial
 gives contribution $3$ to the $x_2$-weight
 From $v_{223}$ the possible values for the $x_2$-weight are even,
 and we never make $6$, again.
 \item[iii)] The monomial contains $v_{012}^6$ and does not contain $v_{000}$.
This monomial gives contribution $6$ to the $x_0$-weight , and the same contribution
is given to the $x_1$-weight and to the $x_2$-weight.
 Hence the only other possible monomial
 that we are allowed to use is $v_{334}$, which gives a $x_3$-weight  doubled with respect to the $x_4$-weight,
which is a contradiction.
   \end{enumerate} 
This contradiction proves our claim. Nevertheless, if our polynomial
is reducible, also its factors have to be homogeneous and $SL(V)$-invariant, and the monomial in the statement
should split into two factors of degree $5$ and $10$, against the claim.\qedd 

 {\it Proof of \thmref{main}} Let $\phi\in \sigma_7(\P {{}}(V),\O(3))$.
By the definition of higher secant variety, $\phi$ is
 in the closure of elements
 which can be written as $\sum_{i=1}^7\phi_i$ with $\phi_i\in (\P {{}}(V),\O(3))$.
 From \lemref{im2}
 it follows that  $$\textrm{rk\ }B_{\phi}\le \textrm{rk\ }B_{\sum_{i=1}^7\phi_i}=
 \textrm{rk\ }\sum_{i=1}^7 B_{\phi_i}\le \sum_{i=1}^7\textrm{rk\ }B_{\phi_i}= 6\cdot 7 = 42$$
 Hence $\det(B_{\phi})$ has to vanish on $\sigma_7(\P {{}}(V),\O(3))$.
 
 We order the monomial basis of $S^3V$ with the lexicographical ordered induced by
 $x_0<x_1<x_2<x_3<x_4$. We order also the basis of
 $\wedge^2V\otimes\wedge^4V$ with the lexicographical order. There are $50$ terms,
 beginning with $$(e_0\wedge e_1)\otimes(e_0\wedge e_1\wedge e_2\wedge e_3),
 (e_0\wedge e_1)\otimes(e_0\wedge e_1\wedge e_2\wedge e_4),\ldots $$
and ending with
 $$\ldots, (e_3\wedge e_4)\otimes(e_1\wedge e_2\wedge e_3\wedge e_4)$$
These $50$ terms are divided into $10$ blocks, depending on the first factor $e_s\wedge e_t$.
The matrix of $B'_{\phi}$,  with respect to this basis, is
 a $50\times 50$ symmetric matrix with linear monomial entries from
 $v_{ijk}$.

We describe this matrix in block form.
For $i=0,\ldots ,4$ let $A_i$ be the $5\times 5$ symmetric matrix which at the entry
$(5-s,5-t)$ has $(-1)^{s+t}v_{ist}$, corresponding to the monomial $x_ix_sx_t$.
For example

$$A_4=\left[\begin{array}{rrrrr}
v_{444}&-v_{344}&v_{244}&-v_{144}&v_{044}\\
-v_{344}&v_{334}&-v_{234}&v_{134}&-v_{034}\\
v_{334}&-v_{234}&v_{224}&-v_{124}&v_{024}\\
-v_{144}&v_{134}&-v_{124}&v_{114}&-v_{014}\\
v_{044}&-v_{034}&v_{024}&-v_{014}&v_{004}\\
\end{array}\right]$$
Then the matrix of $B'_{\phi}$ has the following block form
$$\left[\begin{array}{rrrrrrrrrr}
&&&&&&&A_4&-A_3&A_2\\
&&&&&-A_4&A_3&&&-A_1\\
&&&&A_4&&-A_2&&A_1\\
&&&&-A_3&A_2&&-A_1\\
&&A_4&-A_3&&&&&&A_0\\
&-A_4&&A_2&&&&&-A_0\\
&A_3&-A_2&&&&&A_0\\
A_4&&&-A_1&&&A_0\\
-A_3&&A_1&&&-A_0\\
A_2&-A_1&&&A_0\\
\end{array}\right]$$
 
 Among the $50$ basis elements, there are $30$ tensors  
$(e_s\wedge e_t)\otimes(e_i\wedge e_j\wedge e_k\wedge e_l)$
 such that $\{s,t\}\subseteq \{i,j,k,l\}$. The other $20$ elements are divided into
$5$ groups, depending on the single index $\{s,t\}\cap \{i,j,k,l\}$. The contraction $c$ maps
the first group of $30$ elements into $30$ independent elements of $\Gamma^{2,2,1,1}V$,
and each group of $4$ elements has the image through $c$ of dimension $3$ in 
$\Gamma^{2,2,1,1}V$, indeed the images of the $4$ elements satisfy a linear relation with $\pm 1$ coefficients.

It follows that the matrix of $B_{\phi}$ can be obtained 
from the matrix of $B'_{\phi}$ by deleting five rows, one for each
of the above groups,  and the corresponding five columns.
We can delete, for example, the columns and the rows corresponding to
 
$(e_0\wedge e_1)\otimes(e_1\wedge e_2\wedge e_3\wedge e_4)$,
 $(e_0\wedge e_2)\otimes(e_1\wedge e_2\wedge e_3\wedge e_4)$,
 $(e_0\wedge e_3)\otimes(e_1\wedge e_2\wedge e_3\wedge e_4)$,
 
$(e_0\wedge e_4)\otimes(e_0\wedge e_1\wedge e_2\wedge e_3)$,
 $(e_0\wedge e_4)\otimes(e_1\wedge e_2\wedge e_3\wedge e_4)$

which have respectively number $5, 10, 15, 16, 20$.
Note that in the resulting matrix for $B_{\phi}$, all entries are monomials
in $v_{ijk}$ with coefficient $\pm 1$.
 
 In order to show that for general $\phi$ the morphism $B_{\phi}$ is invertible,
 the simplest way is to look at the monomial $\left(v_{001}v_{022}v_{113}v_{244}v_{334}\right)^9$
 which appears with nonzero coefficient
 in the expression of $\det B_{\phi}$.
 We prefer instead to use the monomial appearing in the statement of \lemref{monomial}, which allows to prove 
  the stronger statement that 
 $\det B_{\phi}$ is the cube of an irreducible polynomial.
 Indeed, by substituting $0$ to all the variables different from
 $v_{000}, v_{012}, v_{111}, v_{223}, v_{334}, v_{144}$,
 we get by an explicit computation that the determinant is equal to
 $$-2\left(v_{000}^2v_{012}^3v_{111}v_{223}^3v_{334}^3v_{144}^3\right)^3$$
 Hence for general $\phi$ we have $\textrm{rk }B_{\phi}=45$.
 Note that this gives an alternative proof of the fact that
 $\sigma_7(\P {{}}(V),\O(3))$ has codimension bigger than zero,
 and it has to appear in the Alexander-Hirschowitz classification.
 It follows that on the points of $\sigma_7(\P {{}}(V),\O(3))$
 the rank of $\textrm{rk }B_{\phi}$ drops at least by three, so that
 $\sigma_7(\P {{}}(V),\O(3))$ is contained in the singular locus of
 $\det B_{\phi}$, and in particular $\det B_{\phi}$ has to vanish
 with multiplicity $\ge 3$ on $\sigma_7(\P {{}}(V),\O(3))$. 
 It is known that
 $\sigma_7(\P {{}}(V),\O(3))$ is a hypersurface (see \cite{CH}), hence its
  equation $P$ has to be a factor of multiplicity $\ge 3$ of
 $\det B_{\phi}$. Since every $SL(V)$-invariant polynomial has degree $5k$, the possible
 values for the degree of $P$ are $5$, $10$ or $15$. Look at the monomials in $P$
 containing some among the variables $v_{000}, v_{012}, v_{111}, v_{223}, v_{334}, v_{144}$,
these monomials have to exist, due to the explicit computation performed before.
 If the degree of $P$ is $\le 10$, then there exists a $SL(V)$-invariant polynomial of degre $10$
 with a monomial containing the above variables, but this contradicts the claim
proved along the proof of
 the \lemref{monomial}. It follows that $\deg P=\deg\sigma_7(\P {{}}(V),\O(3))=15$
 and $P^3$ divides  $\det B_{\phi}$, looking again at our explicit computation
we see that we can arrange the scalar multiples in order that $P$ is defined over the rational numbers
(as all the $SL(V)$-invariants) and the equation $2P(\phi)^3=\det B_{\phi}$ holds.
The \lemref{monomial} shows that $P$ is irreducible.
\qedd

\begin{rem0}\label{reich}
The results obtained by Reichstein with his algorithm developed in \cite{Re}
can be verified with the \thmref{main}. For example when $w$ is like in the example 1
at page 48 of \cite{Re}, a computer check shows that $\textrm{rk\ }(B_w)=42$, confirming that
$w\in\sigma_7(\P {{}}(V),\O(3))$, while
when $w$ is like in the example 2
at page 57 of \cite{Re} then $\textrm{rk\ }(B_w)=45$, so that $w\notin\sigma_7(\P {{}}(V),\O(3))$.

The simplest example of a cubic which is not the sum of seven cubes is probably
$$\phi=x_0^2x_1+x_0x_2^2+x_1^2x_3+x_2x_4^2+x_3^2x_4$$
where $\det(B_{\phi})=-2$, which can be checked even without a computer,
but with a good amount of patience. The polynomial $\phi$ defines a smooth cubic $3$-fold.
\end{rem0}
 
\section{The catalecticant invariant for Clebsch quartics}

Let $U$ be any vector space of dimension $n+1$.

Every quartic $f\in S^4U$ induces the contraction $C_f\colon S^2U^{\vee}\to S^2U$.
Clebsch realized in 1861 that if $f\in (\P n,\O(4))$ then $rk A_f=1$ .
Indeed, with the notations of the previous sections,
$$C_{v^4}(u_1u_2)=24u_1(v)u_2(v)v^2$$
is always a scalar multiple of $v^2$.
Clebsch worked in the case $n=2$ but the same result holds for every $n$. 
If $f\in \sigma_k(\P n,\O(4))$, we get that $C_f$ is the limit of
a sum of $k$ matrices of rank one, then $rk C_f\le k$.
The quartic $f$ is called a Clebsch quartic if and only if $\det C_f=0$,
and this equation gives the catalecticant invariant (see \cite{IK} or \cite{DK}).
A matrix description is the following. Let $D_i$ 
for $i=1,\ldots ,\binom{n+2}{2}$ be a basis of differential operators of second order on $U$.
Then $\det\left(D_iD_j f\right)$ is the catalecticant invariant.

If $n=2$, we write
$$f=f_{0000}x_0^4+4f_{0001}x_0^3x_1+6f_{0011}x_0^2x_1^2+\ldots +12f_{0012}x_0^2x_1x_2+\ldots +f_{2222}x_2^4$$

Then the well known expression for the degree $6$ equation of $\sigma_5({\bf P}^2,\O(4))$
is the following  (we choosed the basis $\partial_{00}, \partial_{01}, \partial_{11},
 \partial_{02}, \partial_{12}, \partial_{22}$) 
$$
\det\left[\begin{array}{cccccc}
f_{0000}&f_{0001}&f_{0011}&f_{0002}&f_{0012}&f_{0022}\\
f_{0001}&f_{0011}&f_{0111}&f_{0012}&f_{0112}&f_{0122}\\
f_{0011}&f_{0111}&f_{1111}&f_{0112}&f_{1112}&f_{1122}\\
f_{0002}&f_{0012}&f_{0112}&f_{0022}&f_{0122}&f_{0222}\\
f_{0012}&f_{0112}&f_{1112}&f_{0122}&f_{1122}&f_{1222}\\
f_{0022}&f_{0122}&f_{1122}&f_{0222}&f_{1222}&f_{2222}\\
\end{array}\right]=0$$

The above equation gives the necessary condition to express a quartic homogeneous polynomial in $3$ 
variables as the sum of $5$ fourth powers.
Mukai proves in \cite{Mu} that a general plane quartic is a sum of $6$ fourth powers in $\infty^3$ ways,
parametrized by the Fano $3$-fold $V_{22}$.

The Clebsch quartics give a hypersurface of degree $\binom{n+2}{2}$ in the space of all quartics.

It follows that this hypersurface contains the variety of $k$-secants to $(\P n, \O(4))$
for $k=\left[\binom{n+2}{2}-1\right]=\frac{n(n+3)}{2}$, 
and it is equal to this secant variety for $1\le n\le 4$, which turns out
to be defective for $2\le n\le 4$. Indeed it is a hypersurface while it is expected that it fills
 the ambient space.  This explains why this example appears in the Alexander-Hirschowitz classification.

{\sf

Giorgio Ottaviani\\
Dipartimento di Matematica U. Dini, Universit\`a di Firenze\\
viale Morgagni 67/A,  50134 Firenze, Italy\\
ottavian@math.unifi.it
}

\end{document}